\documentclass[preprint,12pt]{elsarticle}

\usepackage{epsfig}

\usepackage{amssymb}
 \usepackage{amsthm}
  \usepackage{amsfonts}
  \usepackage{amsmath}
  \usepackage{url}
  \usepackage{hyperref}

  \newtheorem{definition}{Definition}
  \newtheorem{lemma}{Lemma}
\newtheorem{theorem}{Theorem}
\newtheorem{corollary}{Corollary}

\newcommand{\bsu}{\boldsymbol{u}}

\newcommand{\bw}{\boldsymbol{w}}

\newcommand{\ee}{\mathbf{e}}

\def\hhenri#1{}
\def\hchristiane#1 {}

\newcommand\blfootnote[1]{%
  \begingroup
  \renewcommand\thefootnote{}\footnote{#1}%
  \addtocounter{footnote}{-1}%
  \endgroup
}

\journal{Mathematics and Computers in Simulation}

\begin{document}

\begin{frontmatter}



\title{Implementation of irreducible Sobol' sequences in prime power bases}


\author[hf,cl]{Henri Faure and Christiane Lemieux}

\address[hf]{Institut de Math\'ematiques de Marseille,henri.FAURE@univ-amu.fr}
\address[cl]{University of Waterloo, clemieux@uwaterloo.ca}

\begin{abstract}
We present different implementations for the irreducible Sobol' (IS) sequences introduced in \cite{FaLe16}. For this purpose we retain two strategies: first we use the connection between IS and Niederreiter sequences to provide a very simple implementation requiring no computer search; 
 then we use criteria measuring the equidistribution to search for good parameters.
   Numerical results comparing these IS sequences to known implementations 
   show promise for the proposed approaches. 
\end{abstract}


\end{frontmatter}

\blfootnote{\copyright This manuscript version is made available under the CC-BY-NC-ND 4.0 license \url{http://creativecommons.org/licenses/by-nc-nd/4.0/}. Published version has DOI \url{https://doi.org/10.1016/j.matcom.2018.08.015}.}
\section{Introduction}
\label{sec:intro}

Irreducible Sobol' (IS) sequences  \cite{FaLe16}  generalize the famous $LP_\tau$-sequences of Sobol' based on  primitive polynomials over  $\mathbb{F}_2$, to prime power bases and with irreducible polynomials.
This generalization preserves two key properties of Sobol' sequences: $(0,1)$-sequences for their one-di\-men\-sional projections and an easy-to-implement column-by-column construction.
  Just like for 
  $LP_\tau$-sequences, parameters to initialize  the recursions underlying the construction 
  ---the so-called {\em direction numbers}---must be determined.

The goal of this paper is to present different approaches to select direction numbers for IS-sequences and study the properties of their
  implementations, both in terms of their quality and their performance in various numerical integration experiments. In particular we propose an implementation that exploits an insightful connection between IS and Niederreiter sequences. This implementation has the advantage of getting direction numbers ``for free'', without having to search for them. We also propose more traditional implementations based on computer searches for direction numbers. 

This paper is organized as follows. In Section \ref{sec:back} we present important background facts on IS sequences. 
Previous Sobol' sequences constructions are reviewed in Section \ref{sec:prevSob}. Our own implementations are described in Section \ref{sec:implem}, and their quality is assessed in Section \ref{sec:AssessQual}. Numerical results comparing our implementations
 to Sobol' sequences are given in Section \ref{sec:num}. 

\section{Background on IS sequences}
\label{sec:back}


We assume the reader is familiar with the concept of $(t,s)$-sequences, including the definition of the parameter $t$ and equidistribution properties, and refer the reader to \cite{DiPi10} for more information.
We start by recalling the definitions of Sobol' and Niederreiter sequences in the framework of the digital method introduced by Niederreiter and as presented in \cite{FaLe16}. 

The construction introduced by Sobol' in \cite{rSOB67a} is now widely known as {\em Sobol' sequences}. It is a digital sequence in base 2
 very fast to generate, 
 hence its popularity with practitioners. The generating matrices $C^{(i)}$, $1 \le i \le s$, are constructed columns by columns, using monocyclic operators obtained from primitive polynomials over $\mathbb{F}_2$.
To simplify the presentation, we drop the index $i$ and explain below how to construct 
a generating matrix $C$ based on a primitive polynomial $p(x)$ over $\mathbb{F}_2[x]$.

Let $p(x) = a_ex^e + a_{e-1} x^{e-1} + \ldots +a_1 x+a_0$ be a 
primitive polynomial in $\mathbb{F}_2[x]$ of degree $e \ge 1$. 
The matrix $C=(V_r)_{r \ge 1}$ with columns $V_r=(v_{j,r})_{j\ge 1}$ is defined as follows:
For $r=1,\ldots,e$,  let $d_r$ be an odd number between 1 and $2^r$, and let the first $e$ entries of $V_1,\ldots,V_e$ be defined via
\[
\frac{d_r}{2^r} = \sum_{j=1}^{r} v_{j,r} 2^{-j}
\]
and $v_{j,r}=0$ for $j>r$.

Note that since $d_r$ is odd, $v_{r,r}=1$ for $r=1,\ldots,e$. The $e$ integers $d_1,\ldots,d_e$ are the so-called {\em direction numbers} used to initialize the first $e$ column vectors. 
 The remaining column vectors $V_r$ for $r>e$ are obtained using the following linear recurrence associated with $p(x)$:
\begin{equation}
\label{eq:SobRecBase2}
V_{r+e} = \frac{1}{2^e}a_0V_r + a_0V_r + a_1 V_{r+1}+\ldots +a_{e-1} V_{r+e-1}, \qquad r \ge 1,
\end{equation}
where 
 the $j$th entry of $(1/2^e)V_r$ is given by the $(j-e)$th entry of $V_r$ for $j >e$, while the first $e$ entries are 0.

 It is easy to see from \eqref{eq:SobRecBase2} and the property $v_{r,r}=1$ ($1 \le r \le e$) 
 that $C$ is non-singular upper triangular (NUT) and therefore yields a (0,1)-sequence.
 
We point out that Sobol' uses the term {\em direction numbers} for all vectors $V_r$, $r \ge 1$, while we call {\em direction numbers} only the first $e$ ones giving the first $e$ columns of $C$. Hence, 
 the direction numbers associated with $p(x)$ can be defined as the NUT ($e \times e$) {\em direction matrix} $D=(v_{j,r})_{1\le j \le r \le e}$.

As shown in \cite{rSOB67a} (but using a different terminology), Sobol' sequences are $(t,s)$-sequences in base 2 with $t = \sum_{i=1}^s (e_i-1)$, where $e_i$ is the degree of the primitive polynomial used to construct the $i$th generating matrix.

Next, in 1982, Faure introduced $(0,s)$-sequences in a prime base $b \ge s$ (see \cite[Sec.~2.2]{FaLe16}), now widely known as {\em Faure sequences}.  A few years later, 
{\em Niederreiter sequences} were introduced for  a general base $b$  in \cite[Sect.~4]{rNIE88d}.
Here we assume $b$ is  a prime power.
The construction requires $s$ pairwise co-prime polynomials $p_1(x),\ldots,p_s(x) \in \mathbb{F}_b[x]$ of respective positive degrees $e_i$, and then a series of polynomials $g_{i,j}(x) \in \mathbb{F}_b[x]$ for $i=1,\ldots,s$ and $j \ge 1$ such that
$\gcd(p_i(x),g_{i,j}(x)) =1$ for all $i,j$.  The generating matrices are defined through their rows by first developing the formal Laurent series (where $0 \le k <e_i$ and where $w\le 1$ may depend on $i,j,k$)
\begin{equation}
\label{eq:NiedFormal}
\frac{x^k g_{i,j}(x)}{p_i(x)^j}=\sum_{r=w}^{\infty}a^{(i)}(j,k,r)x^{-r}.
\end{equation}
The matrix entries
 are then defined as
$
c_{j,r}^{(i)} = a^{(i)}(q+1,u,r)
$
for $r \ge 1$, where $q$ and $u$
are defined by $j-1=qe_i+u$ with $0 \le u \le e_i-1$.

It is shown in \cite{rNIE88d} that this construction is a digital $(t,s)$-sequence in base $b$ with $t = \sum_{i=1}^s (e_i-1)$, provided that $\lim_{j \rightarrow \infty}(je_i-\deg(g_{i,j}))=\infty$ for all $1 \le i \le s$. This  formula for $t$ is also valid for the Sobol' sequence, 
 but with primitive polynomials. Here however, 
 the polynomials $p_i(x)$ must be co-prime, and thus typically $p_i(x)$ is taken to be the $i$th element in a list of monic irreducible polynomials over 
$\mathbb{F}_b$ sorted in non-decreasing order of degrees, so as to obtain the best possible $t$.
This implies that the parameter $t$ for Niederreiter sequences in base 2 is smaller than $t$ for Sobol' ones.

\subsection{Example showing the relation between Sobol' and Niederreiter sequences}

Consider the primitive polynomial $p(x)=x^2+x+1$ corresponding to the monocyclic linear operator of order 2: $u_{i+2}+u_{i+1}+u_i$ in \cite{rSOB67a}.

In the framework of Sobol', consider the matrix
with starting direction numbers $(1,3)$, resulting from 
 $V_{i+2}=V_{i+1}+V_i+V_i/4$ on column vectors, see \cite[Section 3.2]{rSOB67a}.
In the framework of Niederreiter, consider the matrix 
 generated row by row by the 
series $x^k/p(x)^j$ ($0 \le k < 2$), see \cite[Section 6]{rNIE88d}.

As seen in Figure  \ref{fig:mat_poly_deg2}, a simple examination of these two matrices shows they are the same after permutation of odd and even rows.  Also, it is easy to check on these two matrices that the recurrence relation of Sobol' applies to the original Niederreiter matrix. 
As mentioned before, the Sobol' matrices are NUT matrices and therefore they generate $(0,1)$-sequences. This is an advantage since there is no ``leading-zeros phenomenon" (see \cite[Section 3.3]{rBRA92a}) for Sobol' sequences. Another advantage for implementation is that there is only one recurrence relation for the whole Sobol' matrix instead of a recurrence relation for each odd row of an original Niederreiter matrix in base 2. Hence the interest of a generalization of our example.

\begin{figure}[h]
\[
\begin{bmatrix}
1 & 1 & 0 & 1 & 1 &0 & 1 & 1 & 0 & \ldots\\
0 & 1 & 1 & 0 & 1 & 1 &0 & 1 & 1 & \ldots \\
0 & 0 & 1 & 0 & 1 & 0 & 0 &0 & 1 & \ldots\\
0 & 0 & 0 & 1 & 1 & 1& 0 & 0 &0  & \ldots\\
0 & 0 & 0 & 0 & 1 & 1 & 1 & 0 & 1 & \ldots\\
\vdots & \vdots &  \vdots &  \vdots &  \vdots &  \vdots &  \vdots &  \vdots &  \vdots & \ddots \\
\end{bmatrix}
\qquad
\begin{bmatrix}
0 & 1 & 1 & 0 & 1 &1 & 1 & 0 & 1 & \ldots\\
1 & 1 & 0 & 1 & 1 & 0 &1 & 1 & 0 & \ldots \\
0 & 0 & 0 & 1 & 0 & 1 & 0 &0 & 0 & \ldots\\
0 & 0 & 1 & 0 & 1 & 0 & 0 & 0 & 1 & \ldots\\
0 & 0 & 0 & 0 & 0 & 1 & 1 & 1 & 0 & \ldots\\
\vdots & \vdots &  \vdots &  \vdots &  \vdots &  \vdots &  \vdots &  \vdots &  \vdots & \ddots \\
\end{bmatrix}
\]
\caption{Generating matrices based on $p(x)=x^2+x+1$ for: Sobol' sequence with starting direction numbers $(1,3)$ (left); Niederreiter sequence based on $g_{i,j}=1$ (right).}
\label{fig:mat_poly_deg2}
\end{figure}

\subsection{Definition and properties of IS sequences}

 We now recall the definition of IS sequences introduced in \cite{FaLe16}.

\begin{definition}\label{IrrSob1}{\rm
Let $p(x)=x^e-a_{e-1}x^{e-1}-\cdots-a_1x-a_0$ be a monic irreducible polynomial of degree $e$ over  $\mathbb{F}_b$, where $b$ is a prime power. Define a generating matrix $C$ associated with $p$ by the linear recurrence relation 
\begin{equation}\label{Slrrel}
V_{r+e}-a_{e-1}V_{r+e-1}-\cdots - a_1V_{r+1}-a_0V_r=\frac{1}{b^e}V_r,
\end{equation}
where $V_r$ ($r \ge 1$) is the $r$th column of $C$, and with $e$ starting direction numbers $d_1, \ldots, d_e$ ($1 \le d_r <b^r$ with $\gcd(d_r,b)=1$) defining an NUT ($e \times e$) direction matrix $D$ for $C$.
Then, according to the general principle of construction, an $s$-dimensional {\em irreducible Sobol' sequence} is obtained with $s$ different monic irreducible polynomials $p_i$ generating $s$ such matrices $C^{(i)}$ (typically, one chooses the $e$ first ones in a list of all monic irreducible polynomials sorted in non-decreasing degree, as is done for Niederreiter sequences). 
Note that when working in a general prime power base $b$, one also needs to choose bijections 
to go back and forth between $\mathbb{F}_b$ and $\mathbb{Z}_b$ so that points with coordinates in $[0,1)$ can then be defined.}
\end{definition}

By construction, the generating matrices of irreducible Sobol' sequences  are NUT matrices, so that their one-dimensional projections are $(0,1)$-se\-quen\-ces. Also, it is worth noting that no truncation is required in their definition (in contrast with other types of low-discrepancy sequences). The following lemma and theorems are taken from \cite{FaLe16}.

\begin{lemma} \label{SNrelppb}{\rm  (Fundamental lemma for prime power base $b$ \cite[Lemma 4.2]{FaLe16})} The matrix of a Niederreiter sequence in prime power base  $b$ generated by the formal Laurent series $x^k/p(x)^j$ (for $0 \le k < e$ and $j \ge 1$), where $p$ is a monic irreducible polynomial over $\mathbb{F}_b$ with $\deg(p)=e$, satisfies the Sobol' recurrence relation (\ref{Slrrel}) associated with $p$ in Definition \ref{IrrSob1}.
\end{lemma}

\begin{theorem} (\cite[Theorem 4.3]{FaLe16})
\label{thm:EquivNiedSobolb}
After re-ordering of the rows to get NUT matrices, Niederreiter sequences in  a prime power base $b$ generated by the formal Laurent series $x^k/p_i(x)^j$, where $p_i$, $1 \le i \le s$ are distinct monic irreducible polynomials, 
are IS-sequences associated with the polynomials $p_i$.
\end{theorem}

\begin{theorem}(\cite[Theorem 5.2]{FaLe16})
\label{thm:Nied=IS-iff-gij1}
The only Niederreiter sequences in a prime power base $b$ that are IS-sequences (after re-ordering of the rows to get NUT matrices) are those based on $g_{i,j}(x) = 1$ for all $i=1,\ldots,s$ and $j \ge 1$.
\end{theorem}


\section{Previous Sobol' sequences constructions}
\label{sec:prevSob}

Before we present different implementations for IS sequences, we first review two relatively recent constructions for Sobol' sequences. Both have been defined up to very large dimensions.

The idea proposed by Joe and Kuo \cite{qJOE08a} to find good direction numbers (DNs) is to introduce a criterion that measures the quality parameter $t$ for several two-dimensional projections of the sequence, and then look for the DNs that optimize this criterion via a component-by-component search. The criterion they used to search DNs for coordinate $j$ is 
\begin{equation}
  \label{eq:KJcritParamt}
    {\cal D}^{(q)}_{JK} = \max_{m_{min} \le m \le m_{max}} \frac{[{T}_j(m,w)]^q}{m-{T}_j(m,w)+1},\qquad q>0,
\end{equation}
where 
    $
    T_j(m,w)=\max_{1 \le k < j} (t(j-k,j;m)\times w^{j-k}),
    $
and where $ t(j-k,j;m)$ is the value of the parameter $t$ for the two-dimensional digital net formed by the  $2^m$  first points of the sequence over the coordinates $\{j-k,j\}$, and where $w \in (0,1]$ is a weight typically chosen to be close to 1.

In addition, in their search they also verify if the so-called {\em Property A} holds \cite{rSOB76a}. More precisely,  up to dimension 1111, they only retain DNs that meet this property before assessing them via \eqref{eq:KJcritParamt}.   We say that a sequence satisfies Property $A$ in dimension $s$ if its first $2^s$ points are $(1,\ldots,1)$-equidistributed. That is, if we split the $s$-dimensional hypercube into $2^s$ congruent hypercubes of side $1/2$, then there is one point into each of the $2^s$ sub-cubes of the partition. It should be noted that this property becomes somewhat meaningless once $s$ reaches values beyond which the corresponding number of points $2^s$ is too large to be representative of the number of function evaluations that would be used in real-life problems \cite{qSOB11a}.

In our numerical comparisons we label this sequence `KJ' and have extracted the DN's from Frances Kuo's website, for a sequence built up to 21201 dimensions.

In addition to Property $A$, Sobol' also introduced Property $A^{'}$ in \cite{rSOB76a}, which means that the first $2^{2s}$ points of the sequence are $(2,\ldots,2)$-equidistributed. That is, if we partition the $s$-dimensional hypercube into $2^{2s}$ subcubes of side 0.25, then we have one point in each sub-cube. 

The second type of construction we consider is the one presented in \cite{qSOB11a} under the name SobolSeq16384 (although in this paper we used a version for a slightly lower dimension---6144 instead of 16384---kindly provided to us by S.~Kucherenko, and refer to it as SobolSeq).
This sequence is designed to satisfy Property $A$ up to 6144 dimensions, and Property $A'$ for all sets of five adjacent dimensions (this is referred to as Property $A_5'$ below (borrowing the notation from \cite{qSOB11a})). We note that the latter means that the first $2^{10}=1024$ points of the sequence are tested to see if they are $(2,2,2,2,2)$-equidistributed, over each projection of indices (dimensions) of the form $\{l,l+1,\ldots,l+4\}$, $l=1,\ldots,6140$.

\section{Implementations of IS sequences}
\label{sec:implem}

We now propose two different approaches for implementing IS sequences.

\subsection{Irreducible Sobol'-Niederreiter (ISN) sequences}

This approach simply exploits the connection discussed in Theorem \ref{thm:Nied=IS-iff-gij1}. That is, we construct an IS sequence by defining its corresponding generating matrices as follows. Assume for the $i$th coordinate, the matrix is based on an irreducible polynomial of degree $e$. Then we take the $e$ first rows of the generating matrix of the Niederreiter sequence based on the same polynomial (and with all polynomials $g_{i,j}$ set to 1) and reorder them so that the $e \times e$ upper left generating matrix is NUT. The rest of the matrix is filled column-by-column as in a regular IS sequence. The major advantage of this approach is that we do not need to search for good DNs, and simply need to decide how to order irreducible polynomials of a given degree. In our experiments, we have ordered polynomials $p(x)$ either (1) in increasing order of their decimal representation, or (2) using their decimal representation but interlacing them with the polynomial $\tilde{p}(x)$ with coefficients $a_{e-i}$ for $x^i$. For instance, if $p(x) = x^4+x^3+1$ then $\tilde{p}(x) = x^4+x+1$.  We refer to these two options as ``decimal order'' and ``alternative order'' and label the corresponding sequences as ISN-dec and ISN-alt, respectively.

Note that the generating matrices based on this construction are such that the upper left $e \times e$ matrix is completely determined by the first row, whose $e$ bits are then copied in a diagonal-wise fashion.  This idea will be used later when we describe a method we used to search for ``good'' DNs.

\subsection{Component-by-component search for direction numbers} 
\label{sec:rechParam}

Here we describe two approaches we used to find  ``good'' DNs for IS-sequences in base 2, using a component-by-component search that has some similarities with the approach used in \cite{qJOE08a}. 

$\bullet $ The first approach  explores the space of IS-sequences in base 2 using either the alternative or decimal order. Because that space is very large, we designed a two-step search method. The {\em first step} is that we screen a number of randomly selected DNs (unless the space is small enough to have them all considered) and only retain those who reach the minimum value for a criterion that assesses Property $A_{k_1}$ and Property $A'_{k_2}$, defined in \eqref{eq:crit-propa-ap} below. From the retained DNs, we select the one that minimizes \eqref{eq:ourcritParamt}, a criterion similar to \eqref{eq:KJcritParamt} but that weighs projections differently. 
More precisely, for the first step and assuming we are looking for DNs for coordinate $j$, where $2\le j \le d$,  we first form the $\ell_{j,1}  \times \ell_{j,1}$ matrix obtained by taking the first $\ell_{j,1}=\min(k_1,j)$ elements on the first row of the $\max(1,j-k_1+1)$th,\ldots,$j$th generating matrices and compare its rank  ${\cal R}_{j,k_1}$ with $\ell_{j,1}$. Using the terminology introduced above, if $\ell_{j,1} - {\cal R}_{j,k_1}=0$ for $j=2,\ldots,d$, then we say Property $A_{k_1}$ is met up to dimension $d$.

We then form a $\ell_{j,2 } \!\times\! \ell_{j,2}$ matrix by taking the first $\ell_{j,2}=2\min(k_2,j)$ elements on the first two rows of the  $\max$(1,$j-k_2+1$)th, \ldots,$j$th generating matrices and compare its rank ${\cal R}_{j,k_2}^{'}$ with $\ell_{j,2}$.

Using a weight $\omega \in [0,1]$, we define the criterion
\begin{equation}
  \label{eq:crit-propa-ap}
  \pi_{k_1,k_2} = \omega (\ell_{j,1} - {\cal R}_{j,k_1})+(1-\omega)(\ell_{j,2}-{\cal R}_{j,k_2}^{'}).
  \end{equation}

Then the criterion for the  {\em second step} is of the form 
\begin{equation}
  \label{eq:ourcritParamt}
    {\cal D}^{(q)} = \max_{m_{min} \le m \le m_{max}} \frac{[\hat{T}_j(\ell_2,m,w)]^q}{m-\hat{T}_j(\ell_2,m,w)+1},\qquad q>0
\end{equation}
where 
    $
    \hat{T}_j(\ell_2, m,w)=\max_{1 \le k \le \ell_2} t(j-k,j;m)\times w^{k},
    $
and  where $w \in (0,1]$ is a weight typically chosen  close to 1. We note that ${\cal D}^{(q)}$ depends not just on the parameter $q$ but also on $m_{min},m_{max},\ell_2$ and $w$. 

        So one of the differences with the criterion \eqref{eq:KJcritParamt} from \cite{qJOE08a} is that we do not necessarily look at all two-dimensional projections over indices $\{l,j\}$ with $1 \le l \le j-1$, but instead focus on a window of size $\ell_2$. Doing so reduces the computational burden for calculating the criterion and puts more focus on projections of nearby coordinates, which are more likely to be important when considering the ANOVA decomposition of the function under study. In the above criterion ${\cal D}^{(q)}$, we restrict ourselves to two-dimensional projections, but could easily generalize to a larger set of projections deemed important. 
We label the corresponding construction as IS-t2A-dec/alt in the next two sections, with `dec' or 'alt' referring to the order used for the polynomials.
      
$\bullet $ The second approach focuses on constructions that generalize the  ISN sequences, by considering all possible $e$-bit strings (starting with a 1) for the first $e$ elements of the first row of the direction matrix, 
and choosing the string that minimizes \eqref{eq:ourcritParamt}. We label the obtained sequence as ISN-t2-dec/alt.

Recall that irreducible Sobol'-Niederreiter (ISN) sequences have their upper left $e \times e$ elements completely determined by the first $e$ elements of the first row. In other words, only the most significant bit of the DNs is needed, as we then obtain the other bits on the following $e-1$ rows by shifting the bits by one to the right on each row. Equation \eqref{eq:mat1row} illustrates the process, showing that only the bits $d_2,\ldots, d_{e}$ need to be chosen. (In an ISN sequence, these bits are prescribed by the expansion of $1/p_j(x) = x^{-e}+d_2 x^{-e-1}+d_3 x^{-e-2}+\ldots$).
  \begin{equation}
    \label{eq:mat1row}
    \begin{bmatrix}
      1 & d_2 & d_3 & \cdots & d_e & \cdots \\
      0 & 1 & d_2 & \ddots & d_{e-1} & \cdots \\
      \vdots & \ddots & \ddots &  & \vdots \\
      0 & 0 & 0 & \cdots & 1 & \cdots \\
      \end{bmatrix}
  \end{equation}
Hence only $e-1$ bits need to be chosen instead of $e(e-1)/2$ when all $e$ DNs must be specified, as in Approach 1 or the KJ or SobolSeq constructions.

In the searches that generated the sequences that will be assessed in the next two sections, we have used the criterion ${\cal D}^{(6)}$ with $m_{min}=10,m_{max}=17,\ell_2=20,w=0.9999$ and for the search algorithm that first filters DNs using Property A and A', we used $\omega=0.5$ and $k_1=8,k_2=9$ in \eqref{eq:crit-propa-ap}. The reason for the latter choice is that $k_1=8$ allows us to consider the sequence over a number of points smaller than what is used in ${\cal D}^{(q)}$ (since $m_{min}=10$) and similarly, with $k_2=9$ we are able to assess a point set of size $2^{18}$, hence larger than the largest ones assessed in ${\cal D}^{(q)}$ (since $m_{max}=17$). 

\section{Assessment of quality}
\label{sec:AssessQual}

As should be clear from the previous section, the criteria typically used to measure the quality of Sobol' sequences are the $t$ parameter and Property $A$ (or $A'$). The criteria discussed in the previous section were designed for the search algorithm to select DNs for each component (or dimension) one at a time. Here instead we want to assess and compare sequences over all coordinates from 1 to $d$, and thus the criteria used here are slightly different than in the previous section.  We first describe measures based on the $t$ parameter, and then on Property $A$ and $A'$.
\subsection{Measures based on the $t$-parameter}

For a given value of $m$, we measure the quality of the corresponding point set of size $2^m$ by computing the $t$ parameter denoted $t({\cal J},m)$ over all projections ${\cal J}$ in a  set of the form
\[
{\cal I}(D,d,\mathbf{w}) := \cup_{s=2}^D {\cal I}_{s,w_s,d}
\]
where ${\cal I}_{s,w_s,d}$ contains  ordered $s$-tuples of the form $(i_1,\ldots,i_s)$ with $i_s-i_1+1 \le w_s$, and $i_s \le d$, and $\mathbf{w}$ is a  $(D-1)$-tuple of integers with $w_j \ge j$ for $2 \le j \le D$. 

We can then compute the {\em frequency} vector given by
\[
F_{D,d,\mathbf{w},m} = (n_{0,D,d,\mathbf{w},m},\ldots,n_{m,D,d,\mathbf{w},m}),
\]
where 
$
n_{l,D,d,\mathbf{w},m} = \sum_{s=2}^D \sum_{{\cal J} \in {\cal I}_{s,w_s,d}} \mathbf{1}_{t({\cal J},m) = l}
$
is the number of times we have recorded a value of $l$ for the value of $t({\cal J},m)$ over all subsets ${\cal J}$ considered. 

From the frequency vector, for each $m$ we can compute an average $t$-value 
\[
\bar{t}_{D,d,\mathbf{w},m} = \frac{1}{P_{D,d,\mathbf{w}}}\sum_{l=0}^{m} l \times  n_{l,D,d,\mathbf{w},m}
\]
where $P_{D,d,\mathbf{w}}$ is the cardinality of ${\cal I}(D,d,\mathbf{w})$.

Of course we can also simply look at the maximum value $T_{D,d,\mathbf{w},m}=\max_{{\cal J} \in {\cal I}(D,d,\mathbf{w})} t({\cal J},m)$ obtained over all projections ${\cal J}$ for each $m_0 \le m \le m_1$, and then compute the overall maximum
\[
\tilde{T}_{D,d,\mathbf{w},m_0,m_1}=\max_{m_0 \le m \le m_1} T_{D,d,\mathbf{w},m}.
\]
Finally, we can also compare $t({\cal J},m)$ with its upper bound $\alpha_{{\cal J}}:= \sum_{j \in {\cal J}} (e_j-1)$, where $e_j$ is the degree of the polynomial used in dimension $j$. A nice feature of this measure is that 
$0 \le t({\cal J},m)/\alpha_{{\cal J}} \le 1$, and thus by being scaled it makes it easier to compare this measure across values of $m$. We can even define an overall measure of the form
\[
\tilde{\tau}_{D,d,\mathbf{w}} = \frac{1}{(m_1-m_0+1)P_{D,d,\mathbf{w}}}\sum_{m}\sum_{{\cal J}} \frac{t({\cal J},m)}{\alpha_{{\cal J}}}.
\]

Tables \ref{tab:2da} and \ref{tab:2db} give some results for criteria based on $D=2$. 
For each choice of $(D,d,\mathbf{w})$, we provide the values of $\bar{t}_{D,d,\mathbf{w},m}$ for several $m$ between 4 and 20 (1st line) and $T_{D,d,\mathbf{w},m}$ (second line), its maximal component $\tilde{T}_{D,d,\mathbf{w},m_0,m_1}$ as well as the value $\tilde{\tau}_{D,d,\mathbf{w}}$.

\begin{table}
{\small
\caption{$D=2$, $d=100$,$w_2=100$, 1st line: $\bar{t}_{D,d,\mathbf{w},m}$; 2nd line: $T_{D,d,\mathbf{w},m}$}
\label{tab:2da}
  \begin{tabular}{|l|lllllllll |}
    \hline
    $(\tilde{T}_{D,d,\bw,4,20},\tilde{\tau}_{D,d,\bw}) $ & 4 &  6 &  8 &  10 &  12 &  14 &  16 &  18 &  20  \\
     \hline
   KJ  & 1.3  &  1.9  &  2.4 &  2.7  &  3.0  &  3.2  &  3.3  &  3.4 &   3.6\\
  (8,  0.187) & 3  &5  &6 &7 &7&8&8&8&8 \\
       \hline
   SobolSeq &   1.3 &  1.9&  2.5&  2.8&  3.1&  3.2&  3.3&  3.5&  3.8\\
  (13, 0.201) & 3&5&7&9&11&11&13 &10&10 \\
        \hline
    ISN-alt &1.4&  1.9&  2.3&  2.6& 2.8&    3.0&  3.2&  3.4&  3.5\\
  (11,0.188) &3 &5 &7  &8 &8 &8 &9 &9  &11  \\
        \hline
    ISN-dec&1.4  &  1.9 &  2.3 &  2.6  &  2.8  &  3.0  &  3.2  &  3.4  &  3.5\\
   (11,0.188) & 3  &5  &7 &8  &8  &8 &9 &9  &11\\
        \hline
    ISN-t2-dec  & 1.4 &  1.9 &  2.3  & 2.7 &  2.9 &  3.0 &  3.2  &  3.5 &    3.6\\
  (11,  0.192) & 3 &5  &7  &9  &10  &10  &10  & 10  &11\\
        \hline
    IS-t2A-dec & 1.3 &  2.0 &    2.4 &  2.7  &  3.0 &  3.1 &  3.2  &  3.5  &  3.7\\
(12,  0.196)  & 3 &5  &7 &9 &11  &11  &10  &11  &11\\
   \hline
    \end{tabular}}
    \end{table}

\begin{table}[h]
{\small
\caption{$D=2$, $d=1000$,$w_2=20$, 1st line: $\bar{t}_{D,d,\mathbf{w},m}$; 2nd line: $T_{D,d,\mathbf{w},m}$}
\label{tab:2db}
 \begin{tabular}{|c| l l l l l l l l l |}
        \hline
        $(\tilde{T}_{D,d,\bw,4,20},\tilde{\tau}_{D,d,\bw}) $ & 4 &  6 &  8 &  10 &  12 &  14 &  16 &  18 &  20  \\
     \hline
        KJ & 1.3 &    2.0 &   2.5 &2.9 &    3.2 &   3.5 & 3.7 &    3.9 &   4.1\\
         (12,0.123) & 3 &5 & 7 &9 &10 &11 &11 &12 & 12 \\
        \hline
    SobolSeq &  1.3 & 2.0 &  2.5 &  2.9 &  3.3 & 3.5 &  3.8 & 3.9 &  4.1\\
  (16,0.130) &3  & 5 &7 &9  &11 &13 &15&15 &16  \\
    \hline
    ISN-alt &   1.6 & 2.1 & 2.5 &2.6 & 2.8 & 3.1 &  3.3 & 3.5 &  3.7 \\
 (12,0.120)  & 3  & 5 & 7  & 9 & 10 & 12 & 11 & 12  & 11\\
    \hline
   ISN-dec & 2.3 &   2.6 &  2.5 &   2.5 & 2.8& 3.0 &    3.3 &   3.5 &  3.7 \\
      (12,0.123) & 3&   5 &    7 &   9 &   10 &   10 &   11 &   12 &    12\\
    \hline
    ISN-t2-dec & 1.4&   1.9 & 2.2 &  2.4 &  2.5 &  2.9 &  3.1 &  3.5 &   3.8 \\
     (9,0.113)   & 3 & 5 & 7 & 6 & 6 & 6 & 7 & 8 & 9\\
    \hline
    IS-t2A-dec & 1.3 &    1.9 &    2.3 &     2.5 &   2.8 &    3.1 & 3.4 &     3.8 &    4.0 \\
     (9,0.119) &3 &    5 &      7 &     6 &       6 &     6&      7 &      8 &    9 \\ 
    \hline
    \end{tabular}}
    \end{table}

From these two tables we see that the simple ISN sequences ISN-alt and ISN-dec  often have the best results for the average $t$ value $\bar{t}_{D,d,\mathbf{w},m}$ and the measure $\tilde{\tau}_{D,d,\mathbf{w}}$, while `KJ' does better for the measures based on the maximum value $\tilde{T}_{D,d,\bw,4,20}$ of $t$. We also note that ISN-t2 and IS-t2A perform quite well in Table 2, where $d=1000$ but with the smaller window size $w_2=20$.

\subsection{Measures based on Property A and Property A'}

Using the notation introduced in Section \ref{sec:rechParam}---more precisely the ranks ${\cal R}_{l,k}$ and ${\cal R}'_{l,k}$--- for $d \ge k$, we define the measure
\[
\Pi_{d,k} = \frac{1}{d-1} \sum_{l=2}^{d} (\min(k,l)-{\cal R}_{l,k}).
\]
This corresponds to the average difference between the maximal rank and the actual rank ${\cal R}_{l,k}$ up to dimension $d$. Hence if a sequence is said to satisfy Property $A_k$ up to dimension $d$, then it means $ \Pi_{d,k} =0$. We also introduce
\[
m_{d,k} = \max_{2 \le l \le d} (\min(k,l)-{\cal R}_{l,k})
\]
which returns the largest difference between a rank and its maximum value over all projections considered.

Similar measures are introduced to study Property $A_k'$:
for $d \ge k$ let 
\[
\Pi_{d,k}' = \frac{1}{d-1} \sum_{l=2}^{d} (2\min(k,l)-{\cal R}_{l,k}^{'})
\qquad \mbox{and }
m_{d,k}' = \max_{2 \le l \le d}(2\min(k,l)-{\cal R}_{l,k}^{'}).
\]

In Table \ref{tab:propa} we provide the values related to Property $A$ (that is, ($\Pi_{d,k},m_{d,k}$)) on the first line for each choice of $(d,k)$ and then those related to Property $A^{'}$ (that is, ($\Pi_{d,k}',m_{d,k}')$) on the second line.

\begin{table}[htb]
{\small
\caption{($\Pi_{d,k},m_{d,k}$) (1st line) and   ($\Pi_{d,k}',m_{d,k}')$ (2nd line) for different pairs $(d,k)$}
\label{tab:propa}
    \begin{tabular}{lllllll}
      \hline
    $(d,k)$    & KJ & SobolSeq &  ISN-alt & ISN-dec  & IS-t2A-dec& ISN-t2-alt\\
    \hline
    (100,10) & (0.60,2)&(0.80,2)&(0.94,3)&(1.92,4)&(0.66,2)& (0.80,2)\\
    &(0.80,2) &(0.93,2)&(0.70,2)&(1.22,3)&(0.60,2)& (0.81,2)\\
    \hline
    (360,10) &(0.74,3) &(0.88,3)&(0.89,3)&(0.99,3)&(0.68,2)&(0.81,3)\\
    &(0.78,2)&(0.92,3)& (0.82,2) &(2.53,6)&(0.69,2)&(0.92,3)\\
    \hline
    (1000,10) &(0.77,3) &(0.88,3)&(1.23,3) &(4.00,6)&(0.68,2)&(0.84,3)\\
    &(0.86,3) & (0.85,3) & (1.05,4)& (4.77,9)&(0.67,2)&(0.85,3)\\
    \hline
    (1000,15) & (0.79,3) & (0.87,3) &(1.62,4)  &(5.85,9) &(0.86,3)&(0.82,3)\\
    & (0.82,3) & (0.88,3)&(0.84,3)  & (3.39,8) &(0.83,3)& (0.85,3)\\
    \hline
    (2000,10) &(0.82,3)  & (0.86,3) &(1.63,4) & (4.88,7)&(0.65,2)& (0.82,3)\\
    &(0.85,3) & (0.85,3) & (1.24,4)& (6.20,10) &(0.66,2)& (0.85,3)\\
    \hline
    (5000,10) & (0.84,3) &(0.85,4)  & (2.28,5)  &(6.11,9)&(0.65,2)& (0.79,3)\\
    &(0.85,4) & (0.84,3) & (1.49,5) &(7.79,12)  &(0.68,2)& (0.85,3)\\
    \hline
    \end{tabular}}
    \end{table}


From Table 3 we see a significant difference between the two ways of ordering the polynomials when considering the ISN sequences, with the decimal order performing worse than the alternative one. The difference is even more striking as the dimension increases. This ordering does not matter as much for ISN-t2 and IS-t2A. 
In summary, for Table 3, we see that the construction IS-t2A-dec seems to generally be the best.



\section{Numerical Integration Results}
\label{sec:num}

\hchristiane{Need to add comments}

\begin{figure}[htb]
    \centering
    \includegraphics[width=0.45\textwidth]{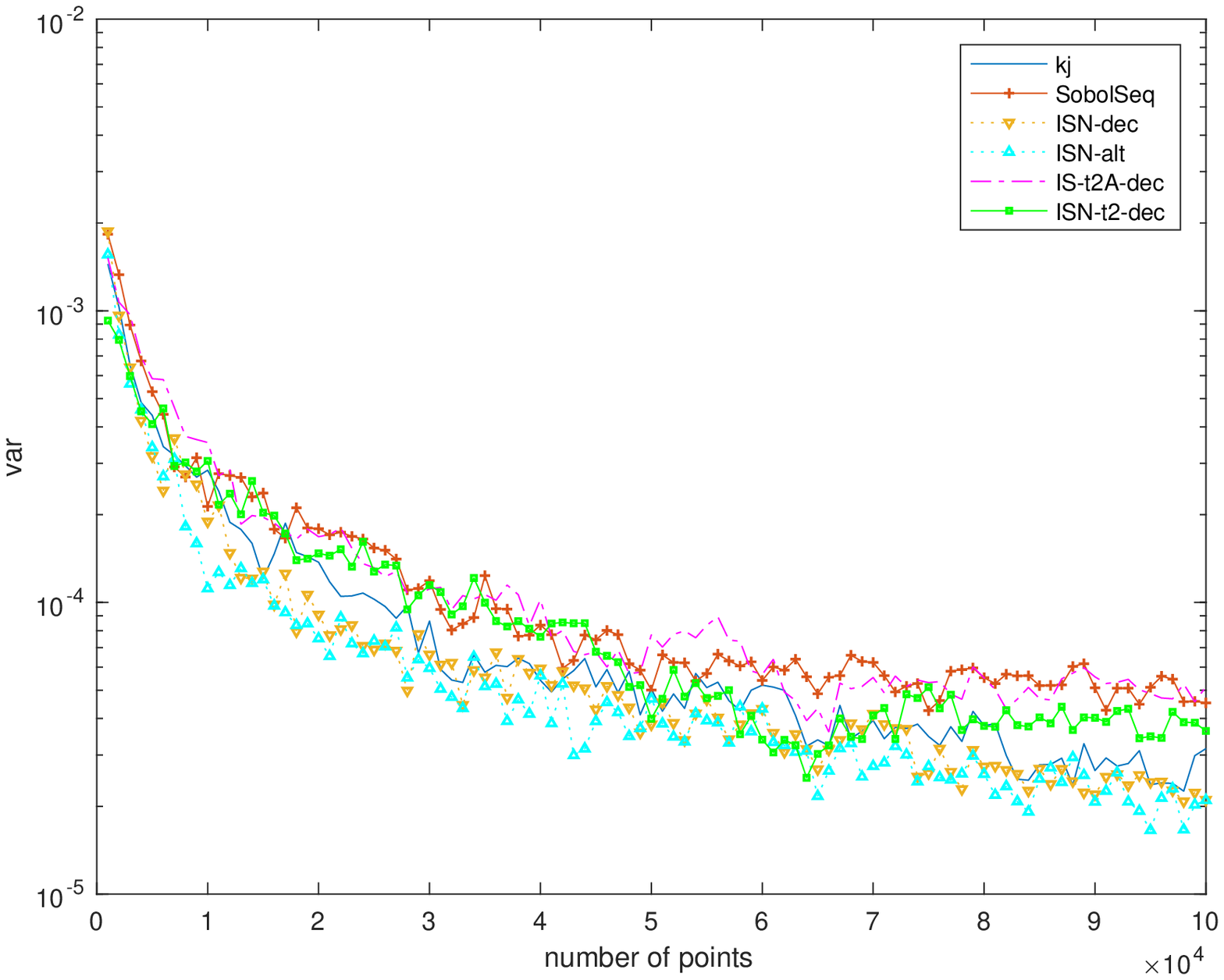}%
    \hfill
    \includegraphics[width=0.45\textwidth]{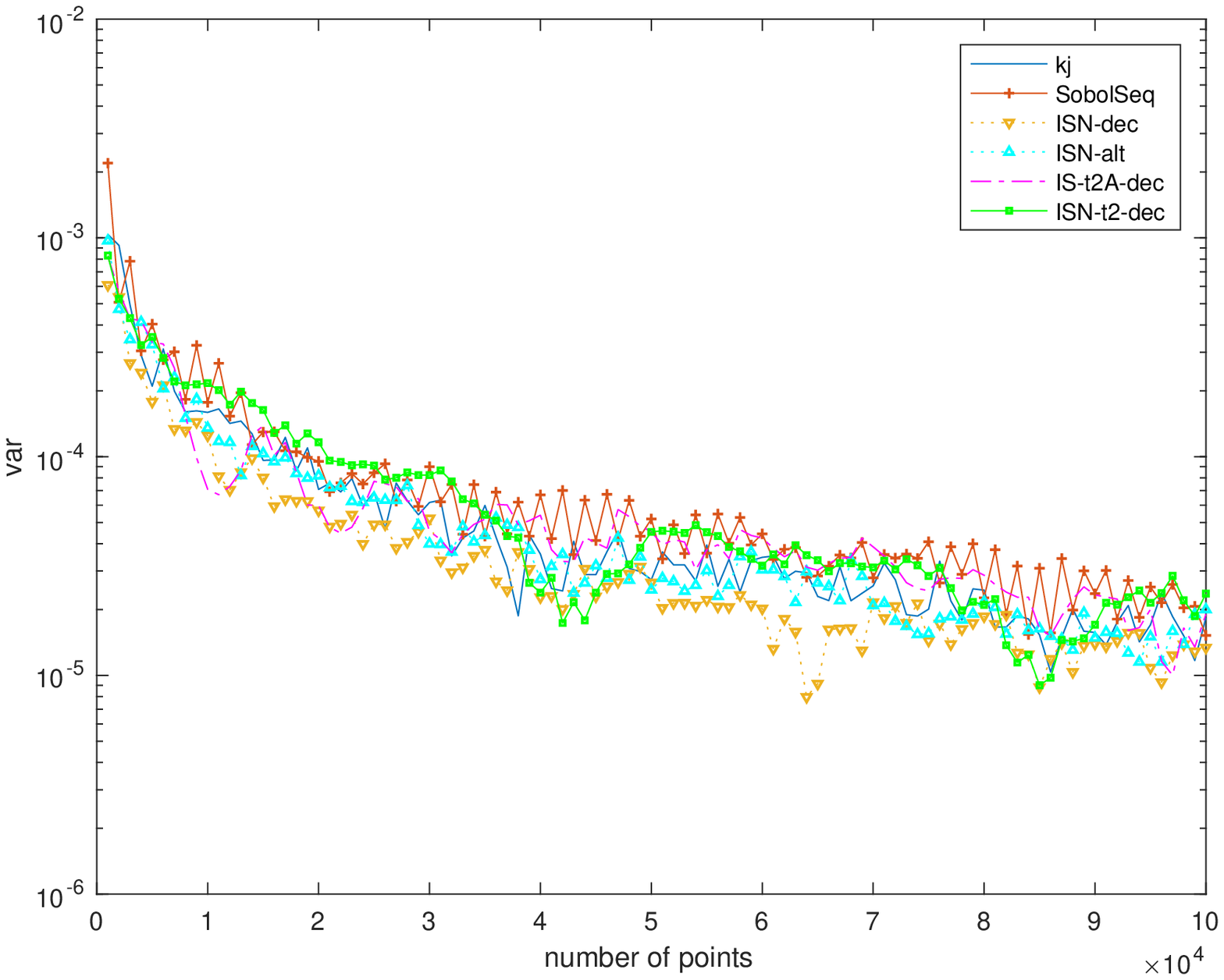}%
    \caption{RMSE for $f_1$: Left: case (i) with $s=1000$; right: case (ii) with $s=20$}
    \label{fig:testfct1}
  \end{figure}

We refer the reader to \cite{qLEM09a} for more information on the functions and examples considered in this section, and  the randomization method used to estimate the variance of the different estimators. First we consider the test function
$
f_1(\bsu) = \prod_{j=1}^s \frac{|4u_j-2|+\alpha_j}{1+\alpha_j}
$
with either the choice (i) $\alpha_j = j$ or (ii) $\alpha_j = s-j+1$. Figure \ref{fig:testfct1} shows the root mean-square error
as a function of the number of points, based on $m=25$ randomizations based on a digital shift. For this problem, the sequences ISN-dec or ISN-alt  
seem to give the smallest error.




Next we consider a problem based on a simple queueing system. Clients arrive according to a Poisson process with arrival rate of 1/minute, and receive service of length that is exponentially distributed with mean 55 seconds. All random variables in this model are assumed to be independent. We simulate the arrival of clients over a fixed period of time $T$ minutes and are interested in $\mathbb{E}(W_{5,T})$, the expected number of clients who will have to wait more than 5 minutes before being served. The problem is thus $2L$-dimensional where $L$ is the number of clients who arrived over $[0,T]$ (a random variable not bounded a priori). We point out that $\mathbb{E}(L) = T$ (with $T$ in minutes). For this  problem, ISN-alt is the best, for both the cases of $T=1000$ and $T=2000$ minutes.
On the right-hand side of Figure \ref{eq:fig4}, we show results from a mortgage-backed security problem often used in the literature, on which ISN-alt, ISN-dec and ISN-t2-dec all do well. 

Overall, based on the results of this section, the simple ISN constructions seem the best. We also note that although ISN-dec does not do well based on the quality measures reported in Table \ref{tab:propa}, these measures consider higher-dimensional projections (up to 15) than what seems to be important in the problems considered in this section, which is why we think this simple construction still does well on those problems.

  \begin{figure}[htbp]
    \centering
    \caption{Queueing problem with $T=1000$ minutes, decimal ordering (left) and alternative ordering (right); shown is the variance of the estimator for $\mathbb{E}(W_{5,T})$}
    \includegraphics[width=0.45\textwidth]{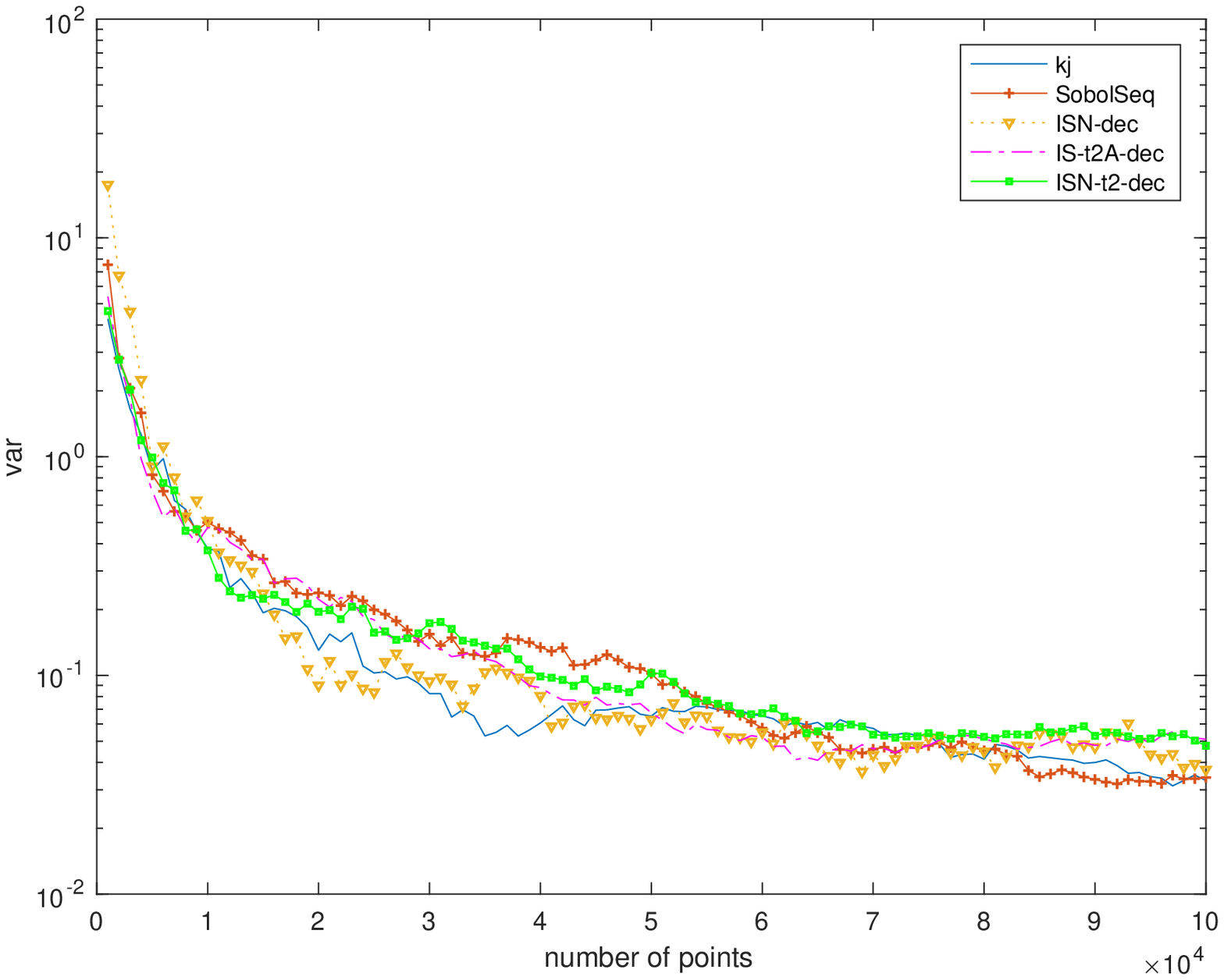}
    \hfill
    \includegraphics[width=0.45\textwidth]{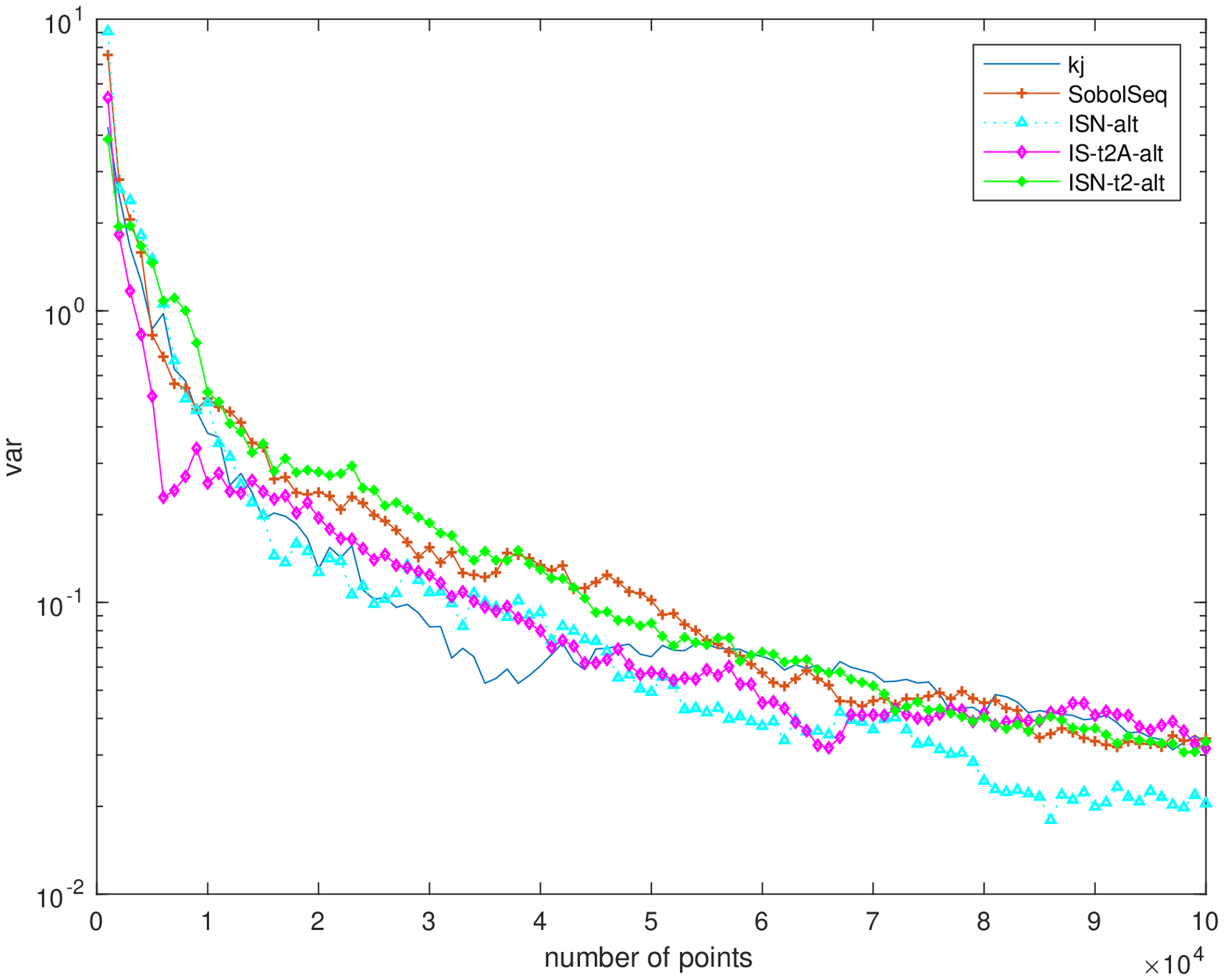}
\end{figure}

\begin{figure}[htb]
\label{eq:fig4}
\centering
\caption{Left: Queueing problem with $T=2000$ minutes; shown is the variance of the estimator for $\mathbb{E}(W_{5,T})$; Right: mortgage-backed security problem in the non-linear setting, shown is the variance of the price at time 0.}
\includegraphics[width=0.45\textwidth]{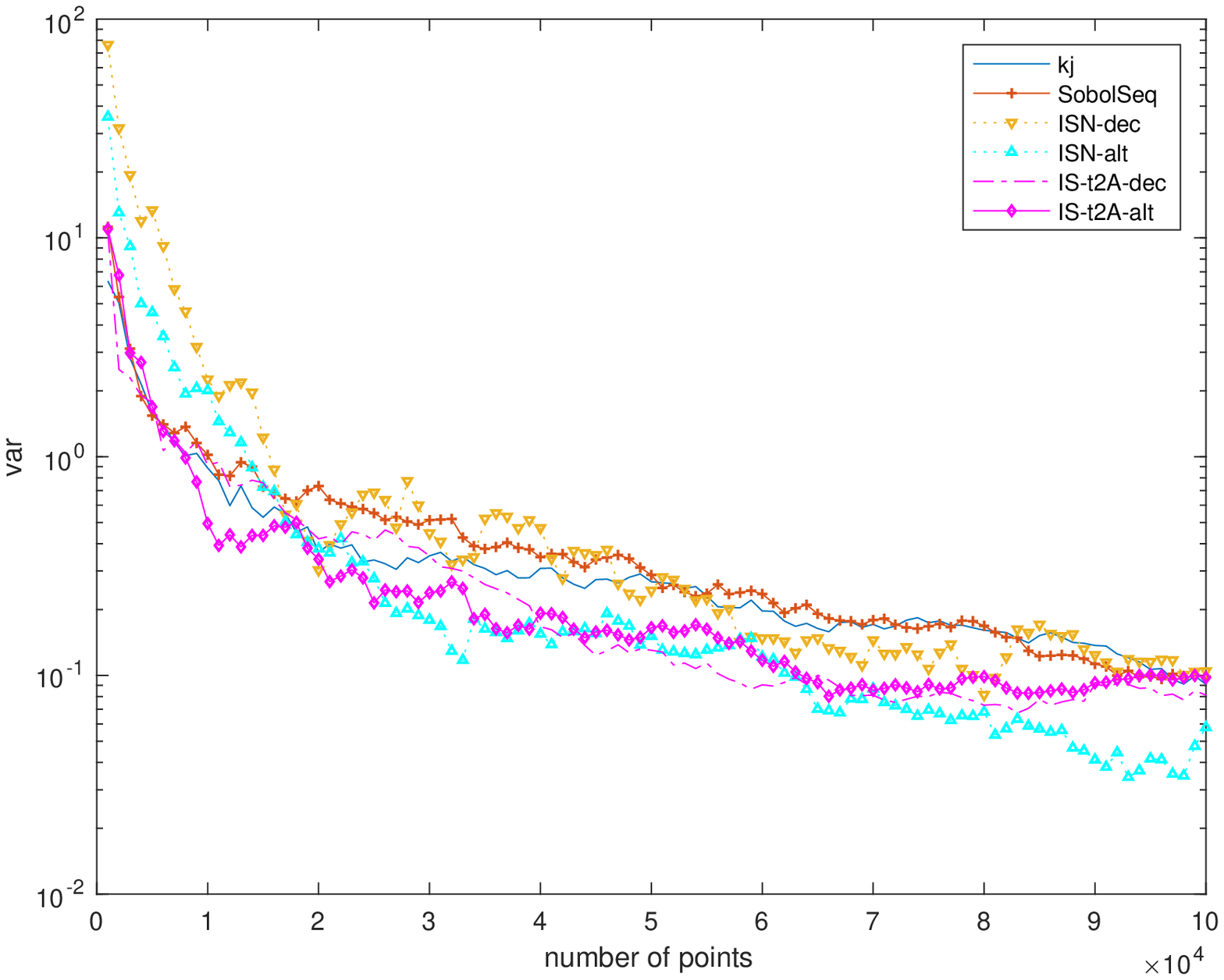}
 \hfill
    \includegraphics[width=0.45\textwidth]{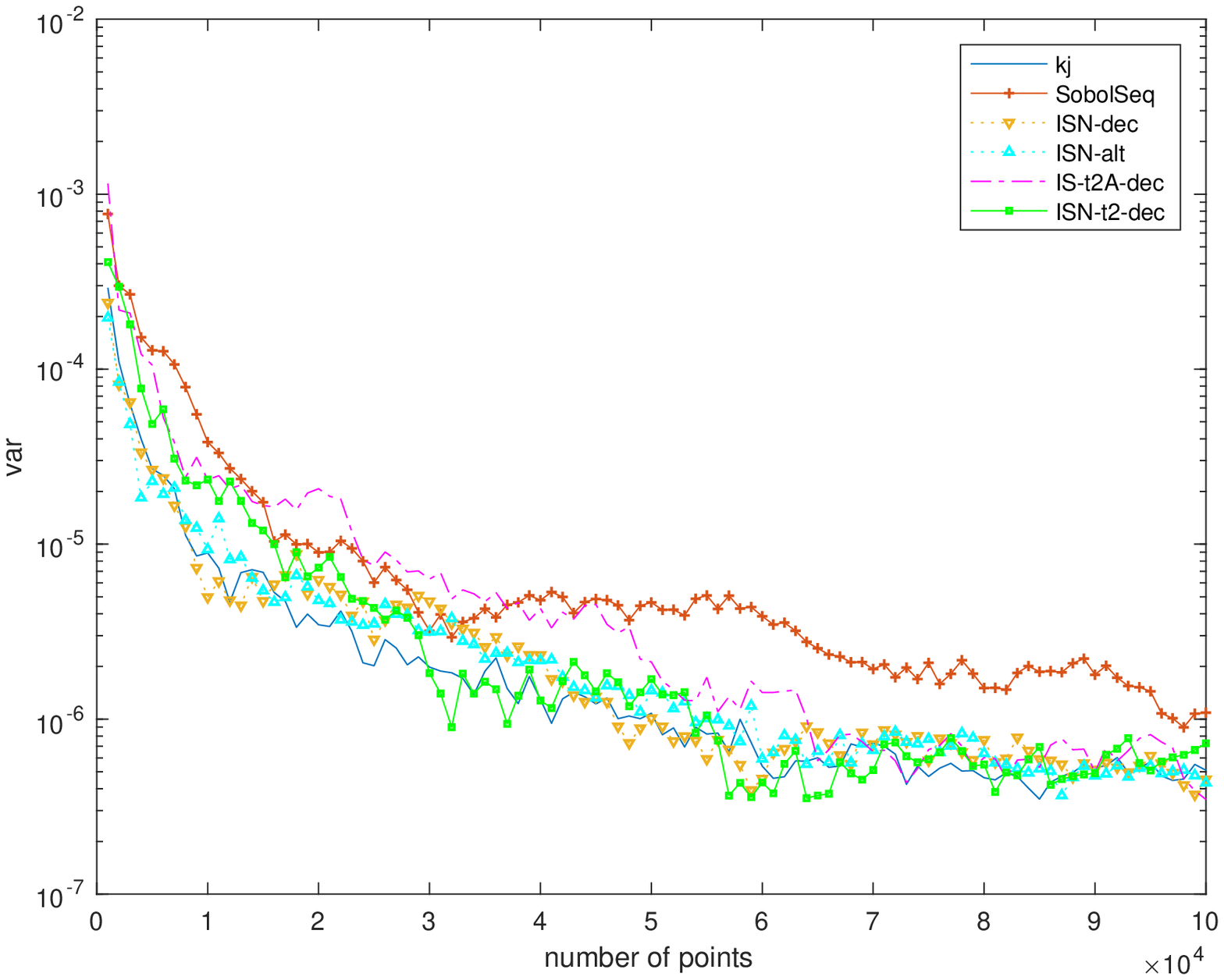}
\end{figure}

\section{Conclusion}

In this paper we have proposed different implementations of IS sequences in base 2. We saw that our  naive irreducible Sobol'-Niederreiter (ISN) implementation---which does not require to search for DN---gives competitive results, even in very high-dimensional problems. This remarkable success of the ISN sequences is very intriguing and a bit of a mystery to us. A theoretical study of this success would be desirable and we plan to pursue it. In particular, it would be interesting to  obtain a family of sequences that include ISN and have the same good results, which could in turn be recommended to users not familiar with the technical background. We also plan to study implementations in bases other than 2.




\section*{Acknowledgments}
We thank the referees for their helpful comments and suggestions, which helped us improve this manuscript. The first author acknowledges the support of NSERC via grant \# 238959.

\end{document}